\documentclass[11pt,a4paper,openany]{article}

\usepackage{amsmath,amsfonts,ams}
\usepackage{latexsym}
\usepackage{amssymb}
\usepackage{amscd}
\usepackage{amsthm}
\usepackage{graphicx}

 \newtheorem{lemma}{Lemma}[section]
 \newtheorem{theorem}{Theorem}[section]
 
 \newtheorem{remark}{Remark}[section]

\begin{document}
\title{Existence theorem and blow-up criterion of the strong solutions to the Magneto-micropolar fluid equations}
\author{{Jia Yuan
$^\dag$ } \\
        {\small The Graduate School of China Academy of Engineering Physics  }\\
        {\small P. O. Box 2101,\ Beijing,\ China,\ 100088 } \\
        {\small (yuanjia930@hotmail.com$^\dag$)}
         \date{}
        }
\maketitle
\begin{abstract}

In this paper we study the magneto-micropolar  fluid equations in
$\R^3$, prove the existence of the strong solution with initial data
in $H^s(\R^3)$ for $s> \frac{3}{2}$, and set up its blow-up
criterion. The tool we mainly use is Littlewood-Paley decomposition,
by which we obtain a Beale-Kato-Majda type blow-up criterion for
smooth solution $(u,\omega,b)$ which relies on the vorticity of
velocity $\nabla\times u$ only.
\end{abstract}

 \vspace{0.2cm}


\section{Introduction}

In this paper,we consider  Magneto-micropolar fluid equations in
$\R^3$.
\begin{eqnarray}\label{1.1}
\begin{cases}
\partial_t u-(\mu+\chi)\Delta u + u \cdot \nabla
u-b\cdot\nabla b+\nabla(p+b^2)-\chi\nabla\times\omega=0, \\
\partial_t \omega-\gamma\Delta\omega-\kappa\nabla
\text{div}\omega+2\chi\omega+u\cdot \nabla\omega-\chi\nabla\times u=0 ,\\
\partial_t b-\nu\Delta b+ u \cdot \nabla b- b \cdot \nabla u=0,\\
 \text{div} u=\text{div} b=0 ,\\
 u(0,x)=u_0(x),\omega(0,x)=\omega_0(x),b(0,x)=b_0(x),
\end{cases}
\end{eqnarray}
where $u(t,x)=(u_1(t,x),u_2(t,x),u_3(t,x)) \in \R^3$ denotes the
velocity of the fluid at a point $x\in\R^3, t\in[0,T)$,
$\omega(t,x)\in\R^3, b(t,x)\in\R^3$ and $p(t,x)\in\R$ denote,
respectively, the micro-rotational velocity, the magnetic field and
the hydrostatic pressure. $\mu,\chi,\kappa,\gamma,\nu$ are positive
 numbers associated to properties of the material: $\mu$ is the kinematic
viscosity, $\chi$ is the vortex viscosity, $\kappa$ and $\gamma$ are
spin viscosities, and $\frac{1}{\nu}$ is the magnetic Reynold.
$u_0,\omega_0,b_0$ are  initial data for the velocity, the angular
velocity and  the magnetic field with properties $\text{div} u_0=0$
and $\text{div} b_0=0$.

There are many earlier results concerning the weak and strong
solvability of magneto-micropolar fluid in  bounded domain
$\Omega\in\R^3$. The corresponding equation is
\begin{eqnarray}\label{1.0}
\begin{cases}
\partial_t u-(\mu+\chi)\Delta u + u \cdot \nabla
u-b\cdot\nabla b+\nabla(p+b^2)-\chi\nabla\times\omega=0, \\
\partial_t \omega-\gamma\Delta\omega-\kappa\nabla
\text{div}\omega+2\chi\omega+u\cdot \nabla\omega-\chi\nabla\times u=0 ,\\
\partial_t b-\nu\Delta b+ u \cdot \nabla b- b \cdot \nabla u=0,\\
 \text{div} u=\text{div} b=0 \qquad\text{in} \Omega, \\
 u(0,x)=u_0(x),\omega(0,x)=\omega_0(x),b(0,x)=b_0(x),\qquad
 x\in\Omega,\\
 u(x,t)=\omega(t,x)=b(t,x)=0, \qquad(t,x)\in[0,T]\times\partial\Omega.
\end{cases}
\end{eqnarray}

 If $b=0$, equation (\ref{1.1})(\ref{1.0}) reduces to
the micropolar fluid system. Micropolar fluid system was first
proposed by Eringe\cite{Er66} in 1966. For the initial
boundary-value problem(\ref{1.0}) with $b=0$, in the year 1977,
Galdi and Rionero\cite{GaRi77} considered the weak solution. Using
linearization and an almost fixed point thereom, in 1988,
Lukaszewicz\cite{Lu88} established the global existence of weak
solutions with sufficiently regular initial data. In 1989, using the
same technique, Lukaszewicz\cite{Lu89} proved the local and global
existence and the uniqueness of the strong solutions. In 2005,
Yamaguchi\cite{Ya05} proved the existence theorem of global in time
solution  for small initial data.

If both $\omega=0$ and $\chi=0$, then the system(\ref{1.1})
reduces to be the magneto-hydrodynamic equations, which has been
studied extensively
in\cite{DuLi72,SeTe83,CaKlSt97,HeXi05,CaMiPrYu06}. Regularity
results can refer to Wu\cite{Wu02, Wu04, Wu06}.

To the full system, Magneto-micropolar fluid equations(\ref{1.0}),
in 1977, Galdi and Rionero\cite{GaRi77} stated the theorem of
existence and uniqueness of strong solutions, but without proof.
Ahmadi and Shaninpoor\cite{AhSh74} studied the stability of
solutions for the system in 1974. By using spectral Galerkin method,
in 1997, Rojas-Medar\cite{RM97} established local existence and
uniqueness of strong solutions. In 1998,Ortega-Torres and
Rojas-Medar\cite{OTRM99} proved global existence of strong solutions
with small initial data. For the weak solution, Rojas-Medar and
Boldrini\cite{RMBo98} established the  local existence in two and
there dimension by  using Galerkin method, and also proved the
uniqueness in 2D case.

But there are few theories about regularity and blow-up criteria of
Magneto-micropolar fluid equations. Some blow-up criterion are
obtained by Yuan\cite{Yu06} in 2006.
 His paper implies that most classical blow-up criteria of smooth solutions to
Navier-Stokes or magneto-hydrodynamic equations also hold for
Magneto-micropolar fluid equations.

For classical MHD equations, an exciting result is that He and
Xin\cite{HeXi05} give a blow-up condition which do not depend on the
magnetic field $b$, which is
$$\int_0^T\|u(t)\|_p^qdt<\infty,\qquad \frac{2}{q}+\frac{3}{p}\leq 1 \qquad 3<p\leq\infty.$$
Zhou\cite{Zh05} gives the regularity criterion dependent only on
$\nabla u$ $$\int_0^T\|\nabla u(t)\|_p^qdt<\infty,\qquad
\frac{2}{q}+\frac{3}{p}\leq 2 \qquad \frac{3}{2}<p<\infty, $$ The
same result s have been extended to Magneto-micropolar fluid
equation by Yuan\cite{Yu06}, of which the condition doesn't rely on
$\omega$ and $b$. We know that for $1<p<\infty$, thanks to the
Biot-Savart law\cite{MaBe02}, $\|\nabla u(t)\|_p$ can be controlled
by $\|\nabla \times u(t)\|_p$, so the regularity criterion of
Zhou\cite{Zh05} can be relaxed by
$$\int_0^T\|\nabla \times u(t)\|_p^qdt<\infty,\qquad
\frac{2}{q}+\frac{3}{p}\leq 2 \qquad \frac{3}{2}<p<\infty,$$ but
this results missed the important marginal case $p=\infty$ which
exactly corresponds to the Beale-Kato-Majda criterion.

While for  $3D$ ideal MHD equations, Caflish, Klapper and
Steele\cite{CaKlSt97} extended the well-known result of
Beale-Kato-Majda\cite{BeKaMa84} for the imcompressible Euler
equations to the 3D ideal MHD equations, that is, if
$$\int_0^T\big(\|\nabla\times u(t)\|_\infty+\|\nabla\times b(t)\|_\infty\big)dt<\infty,$$
then the smooth solution $(u,b)$ can be extended beyond $t=T$. Zhang
and Liu\cite{ZhLi04} extend the condition to
$$\int_0^T\big(\|\nabla\times u(t)\|_{\dot B^0_{\infty,\infty}}+\|\nabla\times b(t)\|_{\dot B^0_{\infty,\infty}}\big)dt<\infty,$$
Cannon,Chen and Miao\cite{CaChMi06} refined  to the following
\begin{eqnarray}\label{1.3}
\lim_{\varepsilon\rightarrow
0}\sup_{j\in\Z}\int_{T-\varepsilon}^T\big(\|\Delta_j(\nabla\times
u)\|_\infty +\|\Delta_j(\nabla\times b)\|_\infty\big)dt=\delta<M,
\end{eqnarray} where $\Delta_j$ is a frequency localization appeared
in Preliminaries.

The aim of our paper, first is using successive approximation method
to obtain the existence of strong solutions in $\R^3$, then using
Fourier frequency localization to set up  blow-up criterion as
(\ref{1.3}) which relying on $\nabla\times u$ only. Our result is
stated as following:

\begin{theorem}(Main theorem)\label{the1}\\
(i) Local existence: Let $s>\frac{3}{2}$, suppose
$(u_0,\omega_0,b_0)\in H^s(\R^3)$ with
$\text{div}u_0=\text{div}b_0=0$, then there exists a positive time
$T(\|(u_0,\omega_0,b_0)\|_{H^s})$ such that a unique solution
$(u,\omega,b)\in C([0,T);H^s)\cap C^1((0,T);H^s)\cap
C((0,T);H^{s+2})$ of
the system (\ref{1.1}) exists.\\
(ii)Blow-up criterion: Suppose that for $s>\frac{3}{2}$,
$(u,\omega,b)\in C([0,T);H^s)\cap C^1((0,T);H^s)\cap
C((0,T);H^{s+2})$ is the smooth solution to equation(\ref{1.1}). If
there exists an absolute constant $M>0$ such that if
\begin{eqnarray}\label{1.2}
\lim_{\varepsilon\rightarrow 0}\sup_{j\in \Z}
\int_{T-\varepsilon}^T\|{\Delta}_j(\nabla\times u)\|_\infty
dt=\delta<M ,
\end{eqnarray}
then $\delta=0,$ and the solution $(u,\omega,b)$ can be extended
past time $t=T$.If
\begin{eqnarray}
\lim_{\varepsilon\rightarrow 0}\sup_{j\in \Z}
\int_{T-\varepsilon}^T\|{\Delta}_j(\nabla\times u)\|_\infty dt\geq
M,
\end{eqnarray}
then the solution blows up at $t=T$ .

\end{theorem}

\section{Preliminaries}
\setcounter{equation}{0}

In this section we  set our notations and recall the
Littlewood-Paley decomposition,
 and review the so called Beinstein estimate and Commutator
 estimate, which are to be used in the proof of our theorem. In what
 follows positive constants will be denoted by $C$ and will change
 from line to line. If necessary, by $C(*,\cdots,*)$ we denote
 constants depending only on the quantities appearing in
 parentheses.

 Let $\mathcal{S}(\R^3)$ be the Schwartz class of rapidly decreasing functions. Given $f\in\mathcal{S}(\R^3)$, the Fourier transform
 of $f$ defined by
 $$\hat{f}(\xi)=(2\pi)^{-\frac{3}{2}}\int_{\R^3}e^{-ix\cdot \xi}f(x)dx$$
We consider $\chi ,\varphi\in \mathcal{S}(\R^3)$ respectively
support in $B=\{\xi\in \R^3,|\xi|\leq\frac{4}{3}\}$ and
$\mathcal{C}=\{\xi\in \R^3,\frac{3}{4}\leq|\xi|\leq\frac{8}{3}\}$
such that
\begin{align*}
&\chi(\xi) +\sum_{j\geq 0}\varphi(2^{-j}\xi)=1  , \quad \forall \xi\in \R^3\\
&\sum_{j\in \Z}\varphi(2^{-j}\xi)=1,\quad \forall \xi\in
\R^3\backslash\{ 0\},
\end{align*}
Setting $\varphi_j=\varphi(2^{-j}\xi)$, then $supp \varphi_j\cap
supp\varphi_j'=\emptyset$ if $|j-j'|\geq2$ and $supp \chi\cap
supp\varphi_j'=\emptyset$ if $|j-j'|\geq1$. Let $h=F^{-1}\varphi$
and $\tilde{h}=F^{-1}\chi$, the dyadic blocks are defined as follows
\begin{align*}
& {\Delta}_j f=\varphi(2^{-j }D)f=2^{3j}\int_{\R^3}h(2^jy)f(x-y)dy
,\\
& {S_j} f=\sum_{k\leq{j-1}}\Delta_k
f=2^{3j}\int_{\R^3}\tilde{h}(2^jy)f(x-y)dy,\quad  j\in \Z.
\end{align*}
Informally, $\Delta_j=S_{j+1}-S_j$ is a frequency projection to the
annulus ${|\xi|\approx2^j}$, while $S_j$ be frequency projection to
the ball ${|\xi|\lesssim2^j}.$ The details of Littlewood-Paley
decomposition can be found in Tribel\cite{Tr83} and
Chemin\cite{Ch98}. Now Besov spaces in $\R^3$ can be defined as
follows:
\begin{align*}
&\dot B_{p,q}^{s}=\Big\{f\in \mathcal{Z}'(\R^3)\big|\quad
\|f\|_{\dot B_{p,q}^{s}}= \Big(\sum_{j\in
\Z}2^{jsq}\|{\Delta}_jf\|_p^q\Big)^\frac1 q<\infty\Big\}
,q\neq\infty
\\
&\dot B_{p,\infty}^{s}=\Big\{f\in
\mathcal{Z}'(\R^3)\big|\quad\|f\|_{\dot B_{p,\infty}^{s}}=\sup_{j\in
\Z}2^{js} \|{\Delta}_j f\|_p<\infty\Big\}
\end{align*}
where $\mathcal{Z}'$ denotes the dual space of
$\mathcal{Z}=\{f\in\mathcal{S} ;D^\alpha\hat{f}(0)=0 ;
\forall\alpha\in \N^n  \quad \text{multi-index}\}$

Now we introduce well-known Bernstein's Lemma and Commutator
estimate, the proof are ommited  here, we can find the details in
Chemin\cite{Ch98}, Chemin and Lerner\cite{ChL95} and Kato and
Ponce\cite{KaPo88}.

\begin{lemma}\label{lem1}(Bernstein's Lemma) \quad Let $1\leq p\leq q\leq\infty$ .
Assume that $f\in L^p$ , then there exist constant $C,C_1,C_2$
independent of $f,j$ such that
\begin{eqnarray}
\sup_{|\alpha|=k}\|\partial^\alpha f\|_q\leq C
2^{jk+3j(\frac{1}{p}-\frac{1}{q})}\|f\|_p  \qquad
supp\hat{f}\subset\{|\xi|\lesssim 2^j\},\\
C_12^{jk}\|f\|_p \leq \sup_{|\alpha|=k}\|\partial^\alpha f\|_p \leq
C_2 2^{jk}\|f\|_p   \quad  supp\hat{f}\subset\{|\xi|\approx 2^j\} .
\end{eqnarray}
\end{lemma}
\begin{remark}
From the above Beinstein estimate, we easily know that in $\R^3$,
for the Reisz transform $R_k(k=1,2,3)$, it has for $\forall 1\leq p
\leq q\leq\infty$
\begin{eqnarray}\label{2.3}
\|R_k\Delta_j u\|_q\leq C2^{3j(\frac{1}{p}-\frac{1}{q})}\|u\|_p.
\end{eqnarray}
 If
suppose vector valued funtion $u$ be divergence free, by Biot Savard
law $\nabla u=(-\Delta)^{-1}\nabla \nabla\times v$ with
$v=\nabla\times u$ and the boundedness of Reisz transform on
$L^p(1<p<\infty)$, we have, there exist constants $C$ independent
$u$ such that
\begin{eqnarray}\label{2.4}
\|\nabla u\|_p\leq C\|v\|_p, \quad  \forall 1<p<\infty.
\end{eqnarray}
If the frequency of $u$ is restricted to annulus
${|\xi|\approx2^j}$, then (\ref{2.3}) implies that
\begin{eqnarray}\label{2.5}
\|\nabla u\|_p\leq C\|v\|_p, \quad  \forall 1\leq p\leq \infty.
\end{eqnarray}
\end{remark}
Now we denote $\Lambda=(I-\Delta)^{\frac{1}{2}}$, which satisfies
$$\widehat{\Lambda f}(\xi)=(1+|\xi|^2)^{\frac{1}{2}}\hat{f}(\xi),$$
 $\Lambda^s(s\in\R)$ can be defined  in the same way
$$\widehat{\Lambda^s f}(\xi)=(1+|\xi|^2)^{\frac{s}{2}}\hat{f}(\xi).$$
Using the above notation, we define the norm of Sobolev space
$W^{s,p}$
$$\|f\|_{W^{s,p}}\triangleq\|\Lambda^s f\|_{L^p},$$
especially  by Fourier transform, $H^s\triangleq W^{s,2}$ can be
defined as
$$H^s\triangleq \big\{f\in\mathcal{S'}(\R^3)\big| \quad\|f\|_{H^s}<\infty \big\},$$
where
$$\|f\|_{H^s}\triangleq\|\Lambda^sf\|_{L^2}=\bigg
(\int_{\R^3}(1+|\xi|^2)^s|\hat{f}(\xi)|^2d\xi\bigg)^{\frac{1}{2}}.$$

\begin{lemma}\label{lem2}(Commutator estimate)\quad Let $1<p<\infty,
s>0$, assume that $f,g\in W^{s,p},$ then there exist constants $C$
independent of $f,g$ such that
\begin{eqnarray}
\big\|[\Lambda^s,f]g\big\|_{L^p}\leq C(\|\nabla
f\|_{L^{p_1}}\|g\|_{W^{s-1,p_2}}+\|f\|_{W^{s,p_3}}\|g\|_{L^{p_4}})
\end{eqnarray}
with $1<p_2,p_3<\infty$ such that
$$\frac{1}{p}=\frac{1}{p_1}+\frac{1}{p_2}=\frac{1}{p_3}+\frac{1}{p_4}.$$
here $[\Lambda^s,f]g=\Lambda^s(fg)-f\Lambda^sg$.
\end{lemma}

\section{Proof of the Theorem \ref{the1} }
\setcounter{equation}{0}

\textbf{Part 1: Local existence}

In order to proof the local existence of equation(\ref{1.1}) with
initial data $(u_0,\omega_0,b_0)\in H^s(\R^3)$ for $s>\frac{3}{2}$,
we construct sequence $(u^{(n+1)},\omega^{(n+1)},b^{(n+1)})$, which
solving the following equations
\begin{eqnarray}
\begin{cases}
\partial_t u^{(n+1)}-(\mu+\chi)\Delta u^{(n+1)}=-u^{(n)} \cdot \nabla
u^{(n)}+b^{(n)}\cdot\nabla b^{(n)}
-\nabla(p^{(n)}+b^{2(n)})\\
\qquad\qquad\qquad\qquad\qquad\qquad\qquad\qquad+\chi\nabla\times\omega^{(n+1)},\\
\partial_t \omega^{(n+1)}-\gamma\Delta\omega^{(n+1)}-\kappa\nabla
\text{div}\omega^{(n+1)}+2\chi\omega^{(n+1)}=-u^{(n)}\cdot \nabla\omega^{(n)}+\chi\nabla\times u^{(n+1)},\\
\partial_t b^{(n+1)}-\nu\Delta b^{(n+1)}=-u^{(n)} \cdot \nabla b^{(n)}+ b^{(n)} \cdot \nabla
u^{(n)},\\
 \text{div} u^{(n+1)}=\text{div} b^{(n+1)}=0 ,\\
 \big(u^{(n+1)},\omega^{(n+1)},b^{(n+1)}\big)(0,x)=S_{n+2}\big(u_0(x),\omega_0(x),b_0(x)\big),
\end{cases}
\end{eqnarray}
for $n=0,1,2,3,\cdots$, where  $b^{2(n)}=(b^{(n)})^2$, and we set
$(u^{(0)},\omega^{(0)},b^{(0)})=(0,0,0)$.

Multiplying the above equation by
$((u^{(n+1)},\omega^{(n+1)},b^{(n+1)}))$ and integrating on time
variable,  denoting $L^2$ inner product by
$\langle\cdot,\cdot\rangle$ ,  we get
\begin{eqnarray}\nonumber
&&\frac{1}{2}\frac{d}{dt}(\|(u^{(n+1)},\omega^{(n+1)},b^{(n+1)})\|_2^2+(\mu+\chi)\|\nabla
u^{(n+1)}\|_2^2+\gamma\|\nabla\omega^{(n+1)}\|_2^2 +\nu\|\nabla
b^{(n+1)}\|_2^2\\ \nonumber
 &&\qquad+\kappa\|div\omega^{(n+1)}\|_2^2+2\chi\|\omega^{(n+1)}\|_2^2\\ \nonumber
&&=-\langle u^{(n)} \cdot \nabla u^{(n)},u^{(n+1)}\rangle+\langle
b^{(n)} \cdot \nabla b^{(n)},u^{(n+1)}\rangle-\langle u^{(n)} \cdot
\nabla \omega^{(n)},\omega^{(n+1)}\rangle\\ \nonumber
 &&\qquad\qquad-\langle u^{(n)} \cdot
\nabla b^{(n)},b^{(n+1)}\rangle+\langle b^{(n)} \cdot \nabla
u^{(n)},b^{(n+1)}\rangle+2\chi\langle \nabla\times
u^{(n+1)},\omega^{(n+1)}\rangle\\
&&=I_1+I_2+I_3+I_4+I_5+I_6,
\end{eqnarray}
where we use $\langle
\nabla\times\omega^{(n+1)},u^{(n+1)}\rangle=\langle \nabla\times
u^{(n+1)},\omega^{(n+1)}\rangle$ .

Using the divergence free condition, the embedding
$H^s\hookrightarrow L^\infty$ and Young Inequality
$$ab\leq\frac{1}{p}a^p+\frac{1}{q}b^q
\qquad\text{where}\qquad \frac{1}{p}+\frac{1}{q}=1,$$ we have
\begin{eqnarray}\nonumber
&&I_1=-\langle u^{(n)} \cdot \nabla
u^{(n)},u^{(n+1)}\rangle\leq\|u^{(n)}\|_2^2\|\nabla
u^{(n+1)}\|_\infty\lesssim \|u^{(n)}\|_2^2\|\nabla
u^{(n+1)}\|_{H^s}\\\nonumber
&&\qquad\qquad\lesssim\frac{\mu}{4}\|\nabla
u^{(n+1)}\|_{H^s}^2+C\|u^{(n)}\|_2^4.
\end{eqnarray}

For the other terms, we use the same technique and get
\begin{eqnarray}\nonumber
&&I_2=\langle b^{(n)} \cdot \nabla
b^{(n)},u^{(n+1)}\rangle\lesssim\frac{\mu}{4}\|\nabla
u^{(n+1)}\|_{H^s}^2+C\|b^{(n)}\|_2^4,\\\nonumber &&I_3=-\langle
u^{(n)} \cdot \nabla
\omega^{(n)},\omega^{(n+1)}\rangle\lesssim\frac{\gamma}{2}\|\nabla
\omega^{(n+1)}\|_{H^s}^2+C\big(\|u^{(n)}\|_2^4+\|\omega^{(n)}\|_2^4\big),\\\nonumber
&&I_4=-\langle u^{(n)} \cdot \nabla
b^{(n)},b^{(n+1)}\rangle\lesssim\frac{\nu}{4}\|\nabla
b^{(n+1)}\|_{H^s}^2+C\big(\|u^{(n)}\|_2^4+\|b^{(n)}\|_2^4\big),\\\nonumber
&&I_5=\langle b^{(n)} \cdot \nabla
u^{(n)},b^{(n+1)}\rangle\lesssim\frac{\nu}{4}\|\nabla
b^{(n+1)}\|_{H^s}^2+C\big(\|b^{(n)}\|_2^4+\|u^{(n)}\|_2^4\big),\\\nonumber
&&I_6=2\chi \langle \nabla\times u^{(n+1)},\omega^{(n+1)}\rangle\leq
\frac{\chi}{2}\|\nabla u^{(n+1)}\|_2^2+2\chi\|\omega^{(n+1)}\|_2^2.
\end{eqnarray}

Summing up the above estimates, we obtain the $L^2$
estimate
\begin{eqnarray}\nonumber
&&\frac{d}{dt}(\|(u^{(n+1)},\omega^{(n+1)},b^{(n+1)})\|_2^2+(\mu+\chi)\|\nabla
u^{(n+1)}\|_2^2+\gamma\|\nabla\omega^{(n+1)}\|_2^2 \\
 &&\quad\quad+\nu\|\nabla b^{(n+1)}\|_2^2
+2\kappa\|div\omega^{(n+1)}\|_2^2\leq
C\|(u^{(n)},\omega^{(n)},b^{(n)})\|_2^4.
\end{eqnarray}

Now let's give the $\dot{H}^s$ estimate. Applying operator
$\Delta_k$ to equation, then multiplying the first three equations
by $(\Delta_k u^{(n+1)},\Delta_k \omega^{(n+1)},\Delta_k
b^{(n+1)})$, introducing notation $\otimes$ as follows
\begin{eqnarray}\nonumber
f\cdot\nabla g=\text{div}(f\otimes g) \quad \text{where} \quad\text{div}(f\otimes
g)^j=\sum_{k=1}^3\partial_k(f^jg^k)=\text{div}(f^jg),
\end{eqnarray}
 we finally
get by using the divergence free condition
\begin{eqnarray}\nonumber\label{3.4}
&&\frac{1}{2}\frac{d}{dt}(\|(\Delta_ku^{(n+1)},\Delta_k\omega^{(n+1)},\Delta_k
b^{(n+1)})\|_2^2+(\mu+\chi)\|\Delta_k\nabla
u^{(n+1)}\|_2^2+\gamma\|\Delta_k\nabla\omega^{(n+1)}\|_2^2 \\
\nonumber && \qquad\qquad+\nu\|\Delta_k\nabla
b^{(n+1)}\|_2^2+\kappa\|\Delta_k \text{div} \omega^{(n+1)}\|_2^2
+2\chi\|\Delta_k\omega^{(n+1)}\|_2^2\\\nonumber
&&=\langle\Delta_k(u^{(n)}\otimes u^{(n)}),\nabla\Delta_k
u^{(n+1)}\rangle-\langle\Delta_k(b^{(n)}\otimes
b^{(n)}),\nabla\Delta_k u^{(n+1)}\rangle\\\nonumber
&&\qquad+\langle\Delta_k(u^{(n)}\otimes \omega^{(n)}),\nabla\Delta_k
\omega^{(n+1)}\rangle+\langle\Delta_k(u^{(n)}\otimes
b^{(n)}),\nabla\Delta_k
b^{(n+1)}\rangle\\
&&\qquad\qquad-\langle\Delta_k(b^{(n)}\otimes
u^{(n)}),\nabla\Delta_k
b^{(n+1)}\rangle-2\chi\langle\Delta_k\omega^{(n+1)},\nabla\times\Delta_k
u^{(n+1)}\rangle,
\end{eqnarray}
where we use the equality
$\langle\Delta_k\omega^{(n+1)},\nabla\times\Delta_k
u^{(n+1)}\rangle=\langle\Delta_k u^{(n+1)},\nabla\times\Delta_k
\omega^{(n+1)}\rangle$.

Multiplying $2^{2ks}$ on both sides of (\ref{3.4}), then summing up
over $k\in\Z$, we  get
\begin{eqnarray}\nonumber
&&\frac{1}{2}\frac{d}{dt}(\|(u^{(n+1)},\omega^{(n+1)},
b^{(n+1)})\|_{\dot{H}^s}^2+(\mu+\chi)\|\nabla
u^{(n+1)}\|_{\dot{H}^s}^2+\gamma\|\nabla\omega^{(n+1)}\|_{\dot{H}^s}^2 \\
\nonumber && \qquad\qquad+\nu\|\nabla
b^{(n+1)}\|_{\dot{H}^s}^2+\kappa\| \text{div}
\omega^{(n+1)}\|_{\dot{H}^s}^2
+2\chi\|\omega^{(n+1)}\|_{\dot{H}^s}^2\\\nonumber
&&\leq\sum_{k\in\Z}2^{2ks}\|\Delta_k(u^{(n)}\otimes
u^{(n)})\|_2\|\Delta_k\nabla
u^{(n+1)}\|_2+\sum_{k\in\Z}2^{2ks}\|\Delta_k(b^{(n)}\otimes
b^{(n)})\|_2\|\Delta_k\nabla u^{(n+1)}\|_2\\\nonumber
&&+\sum_{k\in\Z}2^{2ks}\|\Delta_k(u^{(n)}\otimes
\omega^{(n)})\|_2\|\Delta_k\nabla
\omega^{(n+1)}\|_2+\sum_{k\in\Z}2^{2ks}\|\Delta_k(u^{(n)}\otimes
b^{(n)})\|_2\|\Delta_k\nabla b^{(n+1)}\|_2\\\nonumber
&&+\sum_{k\in\Z}2^{2ks}\|\Delta_k(b^{(n)}\otimes
u^{(n)})\|_2\|\Delta_k\nabla
b^{(n+1)}\|_2+2\chi\sum_{k\in\Z}2^{2ks}\|\Delta_k\nabla\times
u^{(n+1)})\|_2\|\Delta_k\omega^{(n+1)}\|_2\\
&&=II_1+II_2+II_3+II_4+II_5+II_6.
\end{eqnarray}

We use the embedding $H^s\hookrightarrow L^\infty$ along with
H\"{o}lder and Young inequality to get
\begin{eqnarray}\nonumber
II_1\leq\|u^{(n)}u^{(n)}\|_{\dot{H}^s}\|\nabla
u^{(n+1)}\|_{\dot{H}^s}\leq C\|u^{(n)}\|_{L^\infty
}\|u^{(n)}\|_{\dot{H}^s}\|\nabla u^{(n+1)}\|_{\dot{H}^s}\\\nonumber
\lesssim \|u^{(n)}\|_{H^s}^2\|\nabla
u^{(n+1)}\|_{H^s}\lesssim\frac{\mu}{4}\|\nabla
u^{(n+1)}\|_{H^s}^2+C\|u^{(n)}\|_{H^s}^4.
\end{eqnarray}

For the other terms, we use the same technique
\begin{eqnarray}\nonumber
&&II_2\leq\|b^{(n)}b^{(n)}\|_{\dot{H}^s}\|\nabla
u^{(n+1)}\|_{\dot{H}^s}\lesssim\frac{\mu}{4}\|\nabla
u^{(n+1)}\|_{H^s}^2+C\|b^{(n)}\|_{H^s}^4\\\nonumber
&&II_3\leq\|u^{(n)}\omega^{(n)}\|_{\dot{H}^s}\|\nabla
\omega^{(n+1)}\|_{\dot{H}^s}\lesssim\frac{\gamma}{2}\|\nabla
\omega^{(n+1)}\|_{H^s}^2+C\big(\|u^{(n)}\|_{H^s}^4+\|\omega^{(n)}\|_{H^s}^4\big)\\\nonumber
&&II_4+II_5\leq\|b^{(n)}u^{(n)}\|_{\dot{H}^s}\|\nabla
b^{(n+1)}\|_{\dot{H}^s}\lesssim\frac{\nu}{2}\|\nabla
b^{(n+1)}\|_{H^s}^2+C\big(\|u^{(n)}\|_{ H^s}^4+\|b^{(n)}\|_{
H^s}^4\big)\\\nonumber &&II_6\leq2\chi\|\nabla u^{(n+1)}\|_{\dot
H^s}\|\omega^{(n+1)}\|_{\dot H^s}\lesssim\frac{\chi}{2}\|\nabla
u^{(n+1)}\|_{{H}^s}^2+2\chi\|\omega^{(n+1)}\|_{{H}^s}^2
\end{eqnarray}

Taking the sum of $II_1,II_2,II_3,II_4,II_5$ and $II_6$, we obtain
the following estimate
\begin{eqnarray}\nonumber
&&\frac{d}{dt}(\|(u^{(n+1)},\omega^{(n+1)},
b^{(n+1)})\|_{\dot{H}^s}^2+(\mu+\chi)\|\nabla
u^{(n+1)}\|_{\dot{H}^s}^2+\gamma\|\nabla\omega^{(n+1)}\|_{\dot{H}^s}^2 \\
 && \qquad\qquad+\nu\|\nabla
b^{(n+1)}\|_{\dot{H}^s}^2+2\kappa\| \text{div}
\omega^{(n+1)}\|_{\dot{H}^s}^2 \leq C\|(u^{(n)},\omega^{(n)},
b^{(n)})\|_{H^s}^4,
\end{eqnarray}
which along with the $L^2$ estimate, we finally obtain
\begin{eqnarray}\nonumber
&&\frac{d}{dt}(\|(u^{(n+1)},\omega^{(n+1)},
b^{(n+1)})\|_{{H}^s}^2+(\mu+\chi)\|\nabla
u^{(n+1)}\|_{{H}^s}^2+\gamma\|\nabla\omega^{(n+1)}\|_{{H}^s}^2 \\
 && \qquad\qquad+\nu\|\nabla
b^{(n+1)}\|_{{H}^s}^2+2\kappa\| \text{div}
\omega^{(n+1)}\|_{{H}^s}^2 \leq C\|(u^{(n)},\omega^{(n)},
b^{(n)})\|_{H^s}^4.
\end{eqnarray}

Denote
$$E_s^{(n)}(t)=\|(u^{(n)},\omega^{(n)},
b^{(n)})\|_{H^s}^2,$$ then the above inequality can be reduced to be
\begin{eqnarray}\nonumber
&&\frac{d}{dt}E_s^{(n+1)}(t)+(\mu+\chi)\|\nabla
u^{(n+1)}\|_{{H}^s}^2+\gamma\|\nabla\omega^{(n+1)}\|_{{H}^s}^2+\nu\|\nabla
b^{(n+1)}\|_{{H}^s}^2\\
 && \qquad\qquad+2\kappa\| \text{div}
\omega^{(n+1)}\|_{{H}^s}^2  \leq C_1\big(E_s^{(n)}(t)\big)^2.
\end{eqnarray}

Integrating about time variable and taking the supremum  on $[0,T]$,
we have
\begin{eqnarray}\nonumber
&&\sup_{t\in[0,T]}E_s^{(n+1)}(t)+\int_0^T\big((\mu+\chi)\|\nabla
u^{(n+1)}\|_{{H}^s}^2+\gamma\|\nabla\omega^{(n+1)}\|_{{H}^s}^2+\nu\|\nabla
b^{(n+1)}\|_{{H}^s}^2\\\nonumber
 && \qquad\qquad+2\kappa\| \text{div}
\omega^{(n+1)}\|_{{H}^s}^2 \big)dt\\\nonumber && \leq
\|S_{n+2}(u_0,\omega_0,b_0)\|_{H^s}^2+C_1T\big(\sup_{t\in[0,T]}E_s^{(n)}(t)\big)^2\\
&& \leq
C_0\|(u_0,\omega_0,b_0)\|_{H^s}^2+C_1T\big(\sup_{t\in[0,T]}E_s^{(n)}(t)\big)^2.
\end{eqnarray}

By standard induction argument, we find for $\forall n\in\N,
T\in[0,T_0]$, where
\begin{eqnarray}
T_0=\frac{1}{4C_0C_1\|(u_0,\omega_0,b_0)\|_{H^s}^2},
\end{eqnarray}
we can get
\begin{eqnarray}\nonumber
&&\|(u^{(n+1)},\omega^{(n+1)},
b^{(n+1)})\|_{L_T^\infty({H}^s)}+(\mu+\chi)^{\frac{1}{2}}\|\nabla
u^{(n+1)}\|_{L_T^2({H}^s)}+\gamma^{\frac{1}{2}}\|\nabla\omega^{(n+1)}\|_{L_T^2({H}^s)}\\
&&\qquad+\nu^{\frac{1}{2}}\|\nabla
b^{(n+1)}\|_{L_T^2({H}^s)}+(2\kappa)^{\frac{1}{2}}\| \text{div}
\omega^{(n+1)}\|_{L_T^2({H}^s)}\leq2C_0\|(u_0,\omega_0,b_0)\|_{H^s}.
\end{eqnarray}

In the following process, we will show that there exists a positive
time $T_1\leq T$ independent of $n$ such that
$(u^{(n)},\omega^{(n)}, b^{(n)})$ is Cauchy sequence in space
$$\mathcal{X}_{T_1}^{s-1}\triangleq\bigg\{(f,g,h)\in{L_{T_1}^\infty({H}^{s-1})}:\big((\mu+\chi)^{\frac{1}{2}}
\nabla f,\gamma^{\frac{1}{2}}\nabla g, \nu^{\frac{1}{2}}\nabla
h\big)\in {L_{T_1}^2({H}^{s-1})}\bigg\}.$$

Denote $$\delta u^{(n+1)}=u^{(n+1)}-u^{(n)},\delta
\omega^{(n+1)}=\omega^{(n+1)}-\omega^{(n)},\delta
b^{(n+1)}=b^{(n+1)}-b^{(n)}$$
$$\delta p^{(n+1)}=p^{(n+1)}-p^{(n)}, \delta b^{2(n+1)}=(p^{(n+1)})^2-(p^{(n)})^2,$$
which satisfy the following equation
\begin{eqnarray}
\begin{cases}
\partial_t \delta u^{(n+1)}-(\mu+\chi)\Delta \delta u^{(n+1)}=-\delta u^{(n)} \cdot \nabla
u^{(n)}-u^{(n-1)}\cdot\nabla\delta u^{(n)}+\delta b^{(n)}\cdot\nabla b^{(n)}\\
\qquad\qquad\qquad\qquad\qquad+b^{(n-1)}\cdot\nabla\delta b^{(n)}
 -\nabla(\delta p^{(n)}+\delta b^{2(n)})+\chi\nabla\times\delta\omega^{(n+1)}, \\
\partial_t \delta\omega^{(n+1)}-\gamma\Delta\delta\omega^{(n+1)}-\kappa\nabla
\text{div}\delta\omega^{(n+1)}+2\chi\delta\omega^{(n+1)}=-\delta u^{(n)}\cdot \nabla\omega^{(n)}\\
\qquad\qquad\qquad\qquad\qquad-u^{(n-1)}\cdot\nabla\delta\omega^{(n)}+\chi\nabla\times \delta u^{(n+1)},\\
\partial_t \delta b^{(n+1)}-\nu\Delta \delta b^{(n+1)}=-\delta u^{(n)} \cdot \nabla b^{(n)}- u^{(n-1)}\cdot\nabla b^{(n)}
+\delta b^{(n)} \cdot \nabla u^{(n)}\\
\qquad\qquad\qquad\qquad\qquad- b^{(n-1)}\cdot\nabla \delta u^{(n)},\\
 \text{div} \delta u^{(n+1)}=\text{div} \delta b^{(n+1)}=0 ,\\
 \big(\delta u^{(n+1)},\delta \omega^{(n+1)},\delta b^{(n+1)}\big)(0,x)=\Delta_{n+1}\big(u_0(x),\omega_0(x),b_0(x)\big),
\end{cases}
\end{eqnarray}

In the same way, we can get the following estimate
\begin{eqnarray}\nonumber
&&\frac{d}{dt}(\|(\delta u^{(n+1)},\delta\omega^{(n+1)}, \delta
b^{(n+1)})\|_{{H}^s}^2+(\mu+\chi)\|\nabla \delta
u^{(n+1)}\|_{{H}^s}^2+\gamma\|\nabla\delta\omega^{(n+1)}\|_{{H}^s}^2
\\\nonumber
 &&\qquad\qquad+\nu\|\nabla \delta b^{(n+1)}\|_{{H}^s}^2+2\kappa\|
\text{div} \delta\omega^{(n+1)}\|_{{H}^s}^2 \\\nonumber && \leq
C_2\|(\delta u^{(n)},\delta \omega^{(n)},\delta
b^{(n)})\|_2^2\big(\|(u^{(n)},\omega^{(n)},b^{(n)})\|_{H^s}^2+\|(u^{(n-1)},\omega^{(n-1)},
b^{(n-1)})\|_{H^s}^2\big)\\
&&\leq C_3\|(\delta u^{(n)},\delta \omega^{(n)},\delta
b^{(n)})\|_2^2,
\end{eqnarray}
where $C_3=4C_0C_2\|(u_0,\omega_0,b_0)\|_{H^s}$, we uses the
following type of estimates
\begin{eqnarray}\nonumber
 &&\langle\delta u^{(n)}\cdot\nabla u^{(n)},\delta
u^{(n+1)}\rangle\leq\|\delta u^{(n)}\|_2\|
u^{(n)}\|_\infty\|\nabla\delta u^{(n+1)}\|_2\\
&&\qquad\qquad\leq\frac{\nu}{8}\|\nabla\delta
u^{(n+1)}\|_2^2+C\|u^{(n)}\|_{H^s}^2\|\delta
u^{(n)}\|_2^2\\\nonumber
 &&\langle u^{(n-1)}\cdot\nabla\delta u^{(n)},\delta
u^{(n+1)}\rangle\leq\|
u^{(n-1)}\|_\infty\|\delta u^{(n)}\|_2\|\nabla\delta u^{(n+1)}\|_2\\
&&\qquad\qquad\leq\frac{\nu}{8}\|\nabla\delta
u^{(n+1)}\|_2^2+C\|u^{(n-1)}\|_{H^s}^2\|\delta u^{(n)}\|_2^2\\
&&2\chi\langle \nabla\times\delta u^{(n+1)},\delta
\omega^{(n+1)}\rangle\leq\frac{\chi}{2}\|\nabla\delta
u^{(n+1)}\|_2^2+2\chi\|\delta \omega^{(n+1)}\|.
\end{eqnarray}

Integrating over time variable and taking the supremum over $[0,T]$,
denoting $\delta E^{(n)}(t)=\|(\delta u^{(n)},\delta
\omega^{(n)},\delta b^{(n)})\|_2^2$, we obtain
\begin{eqnarray}\nonumber
&&\sup_{t\in[0,T]}\delta
E^{(n+1)}(t)+\int_0^T\big((\mu+\chi)\|\nabla \delta
u^{(n+1)}\|_2^2+\gamma\|\nabla\delta\omega^{(n+1)}\|_2^2+\nu\|\nabla
\delta b^{(n+1)}\|_2^2\\\nonumber
 && \qquad\qquad+2\kappa\| \text{div}
\delta\omega^{(n+1)}\|_2^2 \big)dt\\\nonumber && \leq C_4
2^{-2(n+1)s}\|(u_0,\omega_0,b_0)\|_{H^s}^2+C_3T\delta E^{(n)}(t),
\end{eqnarray}
where we use the fact
$$\|\Delta_{n+1}(u_0,\omega_0,b_0)\|_2^2\leq C_3 2^{-2(n+1)s}\|(u_0,\omega_0,b_0)\|_{H^s}^2.$$

If $C_3 T\leq\frac{1}{2}$, then
\begin{eqnarray}\nonumber
&&\|(\delta u^{(n+1)},\delta\omega^{(n+1)}, \delta
b^{(n+1)})\|_{L_T^\infty(L^2)}+(\mu+\chi)^{\frac{1}{2}}\|\nabla\delta
u^{(n+1)}\|_{L_T^2(L^2)}+\gamma^{\frac{1}{2}}\|\nabla\delta\omega^{(n+1)}\|_{L_T^2(L^2)}\\\nonumber
&&\qquad+\nu^{\frac{1}{2}}\|\nabla\delta
b^{(n+1)}\|_{L_T^2(L^2)}+(2\kappa)^{\frac{1}{2}}\| \text{div}
\delta\omega^{(n+1)}\|_{L_T^2(L^2)}\\
&&\quad\leq2C_4 2^{-(n+1)s}\|(u_0,\omega_0,b_0)\|_{H^s},
\end{eqnarray}
which along with the $H^s$ estimate and the interpolating theorem
$$\|f\|_{H^{s-1}}\leq\|f\|_2^{\frac{1}{s}}\|f\|_{H^s}^{1-\frac{1}{s}},$$
we finally have, by the standard argument, for
$T\leq\min\{T_0,\frac{1}{4C_0}\}$when $n\rightarrow\infty$
\begin{eqnarray}\nonumber
&&\|(\delta u^{(n+1)},\delta\omega^{(n+1)}, \delta
b^{(n+1)})\|_{L_T^\infty(H^{s-1})}+(\mu+\chi)^{\frac{1}{2}}\|\nabla\delta
u^{(n+1)}\|_{L_T^2(H^{s-1})}\\\nonumber
&&+\gamma^{\frac{1}{2}}\|\nabla\delta\omega^{(n+1)}\|_{L_T^2(H^{s-1})}+\nu^{\frac{1}{2}}\|\nabla\delta
b^{(n+1)}\|_{L_T^2(H^{s-1})}+(2\kappa)^{\frac{1}{2}}\| \text{div}
\delta\omega^{(n+1)}\|_{L_T^2(H^{s-1})}\\
&&\quad\leq2C_3^{\frac{1}{s}}C_0^{1-\frac{1}{s}}
2^{-(n+1)}\|(u_0,\omega_0,b_0)\|_{H^s}\rightarrow 0,
\end{eqnarray}
which means $(\delta u^{(n+1)},\delta \omega^{(n+1)},\delta
b^{(n+1)})$ is Cauchy sequence in $\mathcal{X}_{T_1}^{s-1}$, so we
can find the limit $(u,\omega,b)\in\mathcal{X}_{T_1}^{s}$ is a
solution to equation for initial data $(u_0,\omega_0,b_0)\in H^s$,
also the solution satisfies the following estimate
\begin{eqnarray}\nonumber
&&\|(u,\omega,
b)\|_{L_{T_1}^\infty({H}^s)}+(\mu+\chi)^{\frac{1}{2}}\|\nabla
u^{(n+1)}\|_{L_{T_1}^2({H}^s)}+\gamma^{\frac{1}{2}}\|\nabla\omega^{(n+1)}\|_{L_{T_1}^2({H}^s)}\\\nonumber
&&\qquad\qquad+\nu^{\frac{1}{2}}\|\nabla
b^{(n+1)}\|_{L_{T_1}^2({H}^s)}+(2\kappa)^{\frac{1}{2}}\| \text{div}
\omega^{(n+1)}\|_{L_{T_1}^2({H}^s)}\\
&&\leq2C_0\|(u_0,\omega_0,b_0)\|_{H^s}
\end{eqnarray}

This gives the existence of strong solution  of Magneto-micropolar
(\ref{1.1}) in $C([0,T];H^s)$ for $s\geq\frac{3}{2}$. Now let's
prove the uniqueness of the solution.

Suppose $(u,\omega,b),(u',\omega',b')\in L_T^\infty(H^s)$ be two
solutions to equation(\ref{1.1}), let $\delta u=u-u',\delta
\omega=\omega-\omega',\delta b=b-b'$, we deduced that $(\delta
u,\delta \omega,\delta b)$ satisfies the following equation
\begin{eqnarray}
\begin{cases}
\partial_t \delta u-(\mu+\chi)\Delta \delta u=-\delta u \cdot \nabla
u-u'\cdot\nabla\delta u+\delta b\cdot\nabla b+b'\cdot\nabla\delta b
 -\nabla(\delta p+\delta b^2)\\
 \qquad\qquad\qquad\qquad\qquad\qquad\qquad\qquad\qquad\qquad\qquad+\chi\nabla\times\delta\omega,\\
\partial_t \delta\omega-\gamma\Delta\delta\omega-\kappa\nabla
\text{div}\delta\omega+2\chi\delta\omega=-\delta u\cdot \nabla\omega-u'\cdot\nabla\delta\omega+\chi\nabla\times \delta u,\\
\partial_t \delta b-\nu\Delta \delta b=-\delta u \cdot \nabla b- u'\cdot\nabla b
+\delta b \cdot \nabla u- b'\cdot\nabla \delta u,\\
 \text{div} \delta u=\text{div} \delta b=0 ,\\
(\delta u,\delta \omega,\delta b)(0,x)=0.
\end{cases}
\end{eqnarray}

Multiplying the above equation by $(\delta u,\delta \omega,\delta
b)$, then integrating on time variable and using the simple fact

$$\langle u'\cdot\nabla\delta u,\delta u\rangle=\langle u'\cdot\nabla\delta \omega,\delta \omega\rangle
=\langle u'\cdot\nabla\delta b,\delta b\rangle=0$$
$$\langle b'\cdot\nabla\delta b,\delta u\rangle+\langle b'\cdot\nabla\delta u,\delta b\rangle=0,$$
we have
\begin{eqnarray}\nonumber
&&\frac{1}{2}\frac{d}{dt}\|(\delta u,\delta \omega,\delta
b)\|_2^2+(\mu+\chi)\|\nabla\delta u\|_2^2+\gamma\|\nabla\delta
\omega\|_2^2+\nu\|\nabla\delta
b\|_2^2+\kappa\|\text{div}\delta\omega\|_2^2+2\chi\|\delta
\omega\|_2^2\\\nonumber &&=-\langle \delta u\cdot\nabla u,\delta
u\rangle+\langle \delta b\cdot\nabla b,\delta
 u\rangle-\langle \delta u\cdot\nabla \omega,\delta \omega\rangle-\langle \delta u\cdot\nabla b,\delta b\rangle
 +\langle \delta b\cdot\nabla u,\delta b\rangle\\\nonumber
&&\qquad\qquad\qquad\qquad\qquad\qquad\qquad+2\chi\langle
\nabla\times\delta u, \delta
 \omega\rangle\\\nonumber
&&\leq \|\delta u\|_2\|u\|_\infty\|\nabla\delta u\|_2+\|\delta
b\|_2\|b\|_\infty\|\nabla\delta u\|_2+\|\delta
u\|_2\|\omega\|_\infty\|\nabla\delta \omega\|_2\\
&&\qquad\qquad+\|\delta u\|_2\|b\|_\infty\|\nabla\delta
b\|_2+\|\delta b\|_2\|u\|_\infty\|\nabla\delta
b\|_2+2\chi\|\nabla\times\delta u\|_2\|\delta\omega\|_2\\\nonumber
&&\leq \frac{\mu+\chi}{2}\|\nabla\delta
u\|_2^2+\frac{\gamma}{2}\|\nabla\delta
\omega\|_2^2+\frac{\nu}{2}\|\nabla\delta b\|_2^2+2\chi\|\delta
\omega\|_2^2+C\|(u,\omega,b)\|_{H^s}^2\|(\delta
u,\delta\omega,\delta b)\|_2^2,
\end{eqnarray}
that is
\begin{eqnarray}\nonumber
&&\frac{d}{dt}\|(\delta u,\delta \omega,\delta
b)\|_2^2+(\mu+\chi)\|\nabla\delta u\|_2^2+\gamma\|\nabla\delta
\omega\|_2^2+\nu\|\nabla\delta
b\|_2^2+2\kappa\|\text{div}\delta\omega\|_2^2\\\nonumber &&\leq
C\|(u,\omega,b)\|_{H^s}^2\|(\delta u,\delta\omega,\delta b)\|_2^2.
\end{eqnarray}

The $H^s$ estimate imply that
\begin{eqnarray}\nonumber
\|(\delta u,\delta \omega,\delta b)\|_2\leq
2C_0C\|(u_0,\omega_0,b_0)\|_{H^s}T\|(\delta u,\delta\omega,\delta
b)\|_2.
\end{eqnarray}

If $T$ is sufficiently small, we have $\|(\delta u,\delta
\omega,\delta b)\|_2=0$, the proof of local existence is ended up.\\
\\
\textbf{Part 2: Blow-up criterion}

Now we start to proof the second part of Theorem\ref{the1}, to set
up the blow-up criterion. We apply operator $\Lambda^s$  on the two
sides of the equation(\ref{1.1}), multiply $(\Lambda^s u,\Lambda^s
\omega,\Lambda^s b)$ by the resulting equations and integrate the
final form over $\R^3$, and get
\begin{eqnarray}\nonumber
&&\frac{1}{2}\frac{d}{dt}(\|\Lambda^s u\|_2^2+\|\Lambda^s
\omega\|_2^2+\|\Lambda^s b\|_2^2)+(\mu+\chi)\|\nabla\Lambda^s
u\|_2^2+\gamma\|\nabla\Lambda^s \omega\|_2^2 +\nu\|\nabla\Lambda^s
b\|_2^2\\\nonumber
 &&\qquad+\kappa\|div\Lambda^s
\omega\|_2^2+2\chi\|\Lambda^s b\|_2^2\\\nonumber
&&=-\int_{\R^3}\Lambda^s(u\cdot\nabla u)\Lambda^su
dx-\int_{\R^3}\Lambda^s(u\cdot\nabla \omega)\Lambda^s{\omega}
dx-\int_{\R^3}\Lambda^s(u\cdot\nabla b)\Lambda^sb dx\\\nonumber
&&+\int_{\R^3}\Lambda^s(b\cdot\nabla b)\Lambda^su dx
+\int_{\R^3}\Lambda^s(b\cdot\nabla u)\Lambda^sb
dx-2\chi\int_{\R^3}\Lambda^s(\nabla\times u)\Lambda^s \omega dx,
\end{eqnarray}
 where we use the fact
$$\int_{\R^3}\Lambda^s(\nabla\times\omega)\Lambda^s u
dx=\int_{\R^3}\Lambda^s(\nabla\times u)\Lambda^s \omega dx.$$ Now
taking use of the divergence free conditions of $(u,b)$, we have
$$\int_{\R^3}(u\cdot\nabla\Lambda^s u)\Lambda^s udx=\int_{\R^3}(u\cdot\nabla\Lambda^s \omega)\Lambda^s \omega dx
=\int_{\R^3}(u\cdot\nabla\Lambda^s b)\Lambda^s bdx=0,$$ also
applying the following equality
$$\int_{\R^3}(b\cdot\nabla\Lambda^s b)\Lambda^s u
dx+\int_{\R^3}(b\cdot\nabla\Lambda^s u)\Lambda^s b dx=0,$$  we have
\begin{eqnarray}\label{3.2}
\nonumber&&\frac{1}{2}\frac{d}{dt}(\|\Lambda^s u\|_2^2+\|\Lambda^s
\omega\|_2^2+\|\Lambda^s b\|_2^2)+(\mu+\chi)\|\nabla\Lambda^s
u\|_2^2+\gamma\|\nabla\Lambda^s \omega\|_2^2 +\nu\|\nabla\Lambda^s
b\|_2^2\\ \nonumber
 &&\qquad+\kappa\|div\Lambda^s
\omega\|_2^2+2\chi\|\Lambda^s \omega\|_2^2\\ \nonumber
&&=-\int_{\R^3}[\Lambda^s,u]\nabla u\Lambda^su
dx-\int_{\R^3}[\Lambda^s,u]\nabla \omega\Lambda^s\omega dx
-\int_{\R^3}[\Lambda^s,u]\nabla b\Lambda^s b dx \\ \nonumber
&&\qquad+\int_{\R^3}([\Lambda^s,b]\nabla b\Lambda^su+
[\Lambda^s,b]\nabla u\Lambda^sb
)dx-2\chi\int_{\R^3}\Lambda^s(\nabla\times u)\Lambda^s \omega dx.\\
&&=III_1+III_2+III_3+III_4+III_5
\end{eqnarray}
By Lemma\ref{lem2}, H\"{o}lder inequality and Gagliardo-Nirenberg
inequality
\begin{eqnarray}
\|f\|_{W^{s,4}}\leq\|f\|_{W^{s,2}}^{\frac{1}{4}}\|\nabla
f\|_{W^{s,2}}^{\frac{3}{4}},
\end{eqnarray}
 we have
\begin{eqnarray}\label{3.4}
\nonumber&&|III_1|\leq \big\|[\Lambda^s,u]\nabla
u\big\|_{\frac{4}{3}}\|\Lambda^s u\|_4\leq C(\|\nabla u\|_2\|\nabla
u\|_{W^{s-1,4}}+\| u\|_{W^{s,4}}\|\nabla u\|_2)\|u\|_{W^{s,4}}\\
 &&\qquad\leq C\|\nabla
u\|_2\|u\|_{H^s}^{\frac{1}{2}}\|\nabla u\|_{H^s}^{\frac{3}{2}} \leq
C\|\nabla u\|_2^4\| u\|_{H^s}^2+\frac{\mu}{4}\|\nabla u\|_{H^s}^2.
\end{eqnarray}
Using the same technique and the Young inequality
$$ab\leq\frac{1}{p}a^p+\frac{1}{q}b^q ,\qquad\frac{1}{p}+\frac{1}{q}=1,$$
we estimate $III_2,III_3,III_4$ in the same way and get
\begin{eqnarray}\label{3.5}
\nonumber &&|III_2+III_3+III_4|\leq C(\|\nabla u\|_2^4+\|\nabla
\omega\|_2^4+\|\nabla
b\|_2^4)(\|u\|_{H^s}^2+\|\omega\|_{H^s}^2+\|b\|_{H^s}^2)\\
&&\qquad\qquad\qquad\qquad +\frac{\mu}{4}\|\nabla
u\|_{H^s}^2+\frac{\gamma}{2}\|\nabla
\omega\|_{H^s}^2+\frac{\nu}{2}\|\nabla b\|_{H^s}^2.
\end{eqnarray}
Now we estimate the last term
\begin{eqnarray}\label{3.6}
|III_5|\leq 2\chi\|\Lambda^s(\nabla\times
u)\|_2\|\Lambda^s\omega\|_2\leq \frac{\chi}{2}\|\nabla
u\|_{H^s}^2+2\chi\|\omega\|_{H^s}^2.
\end{eqnarray}
Summing up (\ref{3.4})(\ref{3.5})(\ref{3.6}) with (\ref{3.2}), we
get
\begin{eqnarray}\label{3.7}
\nonumber&&\frac{d}{dt}(\|u\|_{H^s}^2+\|\omega\|_{H^s}^2+\|b\|_{H^s}^2)+(\mu+\chi)\|\nabla
u\|_{H^s}^2+\gamma\|\nabla \omega\|_{H^s}^2 +\nu\|\nabla b\|_{H^s}^2
+\kappa\|div \omega\|_{H^s}^2\\
&&\leq
C(\|u\|_{H^1}^4+\|\omega\|_{H^1}^4+\|b\|_{H^1}^4)(\|u\|_{H^s}^2+\|\omega\|_{H^s}^2+\|b\|_{H^s}^2).
\end{eqnarray}
Gronwall inequality gives us
\begin{eqnarray}\label{3.8}
\nonumber&&(\|u\|_{H^s}^2+\|\omega\|_{H^s}^2+\|b\|_{H^s}^2)+\int_0^t\big((\mu+\chi)\|\nabla
u\|_{H^s}^2+\gamma\|\nabla \omega\|_{H^s}^2 +\nu\|\nabla b\|_{H^s}^2
\\\nonumber
&&\qquad\qquad\qquad\qquad\qquad\qquad\qquad\qquad\qquad+\kappa\|div \omega\|_{H^s}^2\big)(t')dt'\\
&&\leq
C(\|u_0\|_{H^s}^2+\|\omega_0\|_{H^s}^2+\|b_0\|_{H^s}^2)\exp\big(t\sup_{t'\in[0,t)}(\|u\|_{H^1}^4+\|\omega\|_{H^1}^4+\|b\|_{H^1}^4)\big).
\end{eqnarray}

Now we go on with the $H^1$ estimates of the solution
$(u,\omega,b)$. Denote $H=\nabla\times
u,I=\nabla\times\omega,J=\nabla\times b$, we
 take curl on both sides of (\ref{1.1}), we get the following
equation
\begin{eqnarray}
\begin{cases}
\partial_t H-(\mu+\chi)\Delta H + u \cdot \nabla
H-H \cdot \nabla u-b\cdot\nabla J+J\cdot \nabla b-\chi\nabla\times I=0, \\
\partial_t I-\gamma\Delta I+2\chi I+u\cdot \nabla I-H\cdot \nabla \omega-\chi\nabla\times H=0 ,\\
\partial_t J-\nu\Delta J+ u \cdot \nabla J- H \cdot \nabla b- b \cdot \nabla H+ J \cdot \nabla u=0
\end{cases}
\end{eqnarray}
which uses the fact $\nabla\times\nabla div \omega=0$.\\
 Multiplying the three equations with $(H,I,J)$ separately,
 integrating over $\R^3$ about the variable $x$, using integrating by parts and the divergence free condition of $u,b$,
 we obtain the fact
 $$\int_{\R^3}(u\cdot\nabla)H\cdot H dx=\int_{\R^3}(u\cdot\nabla)I\cdot I dx=\int_{\R^3}(u\cdot\nabla)J\cdot J dx=0,$$
 $$\int_{\R^3}(b\cdot\nabla)J\cdot H dx+\int_{\R^3}(b\cdot\nabla)H\cdot J dx=0,$$
 $$\int_{\R^3}(\nabla\times H)\cdot I dx=\int_{\R^3}(\nabla\times I)\cdot H dx,$$
 so we finally have
 \begin{eqnarray}\label{3.10}
\nonumber&&\frac{1}{2}\frac{d}{dt}(\|H\|_2^2+\|I
\|_2^2+\|J\|_2^2)+(\mu+\chi)\|\nabla H\|_2^2+\gamma\|\nabla I\|_2^2
+\nu\|\nabla J\|_2^2
 +2\chi\|I\|_2^2\\ \nonumber
&&=\int_{\R^3}(H\cdot\nabla)u\cdot H dx
+\int_{\R^3}(H\cdot\nabla)\omega\cdot I
dx+\int_{\R^3}(J\cdot\nabla)u\cdot J dx\\
\nonumber&&\qquad -\int_{\R^3}(J\cdot\nabla)b\cdot H dx+\int_{\R^3}
(H\cdot\nabla)b\cdot J dx+2\chi\int_{\R^3}(\nabla\times H)\cdot I
dx\\
&&=IV_1+IV_2+IV_3+IV_4+IV_5+IV_6
\end{eqnarray}
Let us first estimate $IV_1$ and $IV_3$,  we use Littlewood-Paley
decomposition to $u$ and dispose it in different frequencies.
\begin{eqnarray}
\nonumber&&IV_1=\sum_{j<-N}\int_{\R^3}(H\cdot\nabla)\Delta_ju\cdot H
dx+\sum_{-N\leq j\leq N}\int_{\R^3}(H\cdot\nabla)\Delta_ju\cdot H
dx\\
\nonumber&&\qquad\qquad\qquad+\sum_{j>N}\int_{\R^3}(H\cdot\nabla)\Delta_ju\cdot
H dx,\\
&&\quad=V_1+V_2+V_3
\end{eqnarray}
For the first term, we have, by H\"{o}lder inequality, Beinstein
inequality and (\ref{2.4}) (\ref{2.5})
\begin{eqnarray}\label{3.12}
|V_1| \leq\|H\|_2^2\sum_{j<-N}\|\nabla\Delta_ju\|_{\infty}\leq
C\|H\|_2^2\sum_{j<-N}2^{\frac{3}{2}j}\|\Delta_jH\|_2 \leq
C2^{-\frac{3}{2}N}\|H\|_2^3 \\
|V_2| \leq\|H\|_2^2\sum_{-N\leq j\leq
N}\|\nabla\Delta_ju\|_{\infty}\leq C \|H\|_2^2\sum_{-N\leq j\leq
N}\|\Delta_jH\|_{\infty}
\end{eqnarray}
for $V_3$, using similar method along with interpolation inequality
$$\|H\|_3\leq C\|H\|_2^{\frac{1}{2}}\|\nabla H\|_2^{\frac{1}{2}},$$
we get
\begin{eqnarray}\label{3.14}
 \nonumber|V_3|
&&\leq\|H\|_3^2\sum_{j>N}\|\nabla\Delta_ju\|_3\leq C
\|H\|_3^2\sum_{j>N}2^\frac{j}{2}\|\Delta_jH\|_ 2\\
\nonumber\qquad&&\leq C\|H\|_3^2\bigg(\sum_{j>N}2^{-\frac{j}{2}\cdot
2}\bigg)^\frac{1}{2}\bigg(\sum_{j>N}2^{j\cdot
2}\|\Delta_jH\|_2^2\bigg)^\frac{1}{2}\\
&&\leq C 2^{-\frac{N}{2}}\|H\|_2\|\nabla H\|_2^2.
\end{eqnarray}
Summing up (\ref{3.12})-(\ref{3.14}), we have
\begin{eqnarray}\label{3.15}
|IV_1|\leq C\big(2^{-\frac{3}{2}N}\|H\|_2^3+\|H\|_2^2\sum_{-N\leq
j\leq N}\|\Delta_jH\|_{\infty}+2^{-\frac{N}{2}}\|H\|_2\|\nabla
H\|_2^2\big).
\end{eqnarray}
$IV_3$ can be treated in the same way, we decompose it as
\begin{eqnarray}
\nonumber&&IV_3=\sum_{j<-N}\int_{\R^3}(J\cdot\nabla)\Delta_ju\cdot J
dx+\sum_{-N\leq j\leq N}\int_{\R^3}(J\cdot\nabla)\Delta_ju\cdot J
dx\\
\nonumber&&\qquad\qquad\qquad+\sum_{j>N}\int_{\R^3}(J\cdot\nabla)\Delta_ju\cdot
J dx,
\end{eqnarray}
then obtain the estimate
\begin{eqnarray}\label{3.16}
|IV_3|\leq
C\big(2^{-\frac{3}{2}N}\|J\|_2^2\|H\|_2+\|J\|_2^2\sum_{-N\leq j\leq
N}\|\Delta_jH\|_{\infty}+2^{-\frac{N}{2}}\|J\|_2\|\nabla
J\|_2\|\nabla H\|_2\big).
\end{eqnarray}
Now we study $IV_2,IV_4,IV_5$, we decompose $H$ by using
Littlewood-Paley theory, that is
\begin{eqnarray}
\nonumber&&IV_2=\sum_{j<-N}\int_{\R^3}(\Delta_jH\cdot\nabla)\omega\cdot
I dx+\sum_{-N\leq j\leq
N}\int_{\R^3}(\Delta_jH\cdot\nabla)\omega\cdot I
dx\\
\nonumber&&\qquad\qquad\qquad+\sum_{j>N}\int_{\R^3}(\Delta_jH\cdot\nabla)\omega\cdot
I dx,
\end{eqnarray}
then
\begin{eqnarray}\label{3.17}
|IV_2|\leq
C\big(2^{-\frac{3}{2}N}\|I\|_2^2\|H\|_2+\|I\|_2^2\sum_{-N\leq j\leq
N}\|\Delta_jH\|_{\infty}+2^{-\frac{N}{2}}\|I\|_2\|\nabla
I\|_2\|\nabla H\|_2\big).
\end{eqnarray}
For $IV_4$ and $II_5$,  similarly we have
\begin{eqnarray}\label{3.18}\nonumber
&&|IV_4|+|IV_5|\leq
C\big(2^{-\frac{3}{2}N}\|J\|_2^2\|H\|_2+\|J\|_2^2\sum_{-N\leq j\leq
N}\|\Delta_jH\|_{\infty}\\
&&\qquad\qquad\qquad\qquad\qquad\qquad\qquad\qquad+2^{-\frac{N}{2}}\|J\|_2\|\nabla
J\|_2\|\nabla H\|_2\big).
\end{eqnarray}
Simply using Young inequality, the last term $IV_6$ can be written
as
\begin{eqnarray}\label{3.19}
|IV_6|\leq 2\chi\|\nabla\times
H\|_2\|I\|_2\leq\frac{\chi}{2}\|\nabla H\|_2^2+2\chi\|I\|_2^2.
\end{eqnarray}
Summing up
(\ref{3.15})(\ref{3.16})(\ref{3.17})(\ref{3.18})(\ref{3.19}) and
taking the sum into (\ref{3.10}), by Young inequality, we  get
\begin{eqnarray}\label{3.20}
\nonumber&&\frac{d}{dt}(\|H\|_2^2+\|I
\|_2^2+\|J\|_2^2)+(2\mu+\chi)\|\nabla
H\|_2^2+2\gamma\|\nabla I\|_2^2 +2\nu\|\nabla J\|_2^2\\
\nonumber&&\leq
C\big(2^{-\frac{3}{2}N}(\|H\|_2^3+\|I\|_2^3+\|J\|_2^3)\big)+\sum_{{-N}\leq
j \leq N}\|\Delta_jH\|_{\infty}(\|H\|_2^2+\|I\|_2^2+\|J\|_2^2)\\
&&\qquad+2^{-\frac{N}{2}}(\|H\|_2+\|I\|_2+\|J\|_2)(\|\nabla
H\|_2^2+\|\nabla I\|_2^2+\|\nabla J\|_2^2)
\end{eqnarray}
If we let $2^{-\frac{N}{2}}(\|H\|_2+\|I\|_2+\|J\|_2)\leq
\min(\mu,\gamma,\nu)$, that is, if we choose
\begin{eqnarray}\label{3.21}
N\geq \bigg[\frac{2}{\log
2}log^+\big(\frac{C}{\min(\mu,\gamma,\nu)}(\|H\|_2+\|I\|_2+\|J\|_2)\big)\bigg]+1,
\end{eqnarray} where $[a]$ stands for the integral parts of $a\in\R$,
 $\log^+(x)=\log(x+e)$, then we have
\begin{eqnarray}\label{3.21}
\nonumber&&\frac{d}{dt}(\|H\|_2^2+\|I
\|_2^2+\|J\|_2^2)+(\mu+\chi)\|\nabla
H\|_2^2+\gamma\|\nabla I\|_2^2 +\nu\|\nabla J\|_2^2\\
&&\qquad\qquad\leq C \sum_{{-N}\leq j \leq
N}\|\Delta_jH\|_{\infty}(\|H\|_2^2+\|I\|_2^2+\|J\|_2^2)+C.
\end{eqnarray}
Gronwall inequality gives us that
\begin{eqnarray}
\nonumber&&\|H\|_2+\|I \|_2+\|J\|_2\leq \exp\big(C\sum_{-N\leq j\leq
N}\int_0^t\|\Delta_jH(t')\|_\infty
dt'\big)\\
&&\qquad\qquad\qquad\qquad\qquad\qquad\big(\sqrt{Ct}+\|H(0)\|_2+\|I(0)\|_2+\|J(0)\|_2\big),
\end{eqnarray}
which implies
\begin{eqnarray}\label{3.24}
\nonumber&&\|H\|_2+\|I\|_2+\|J\|_2\leq
\exp\big(C\log^+(\|H\|_2+\|I\|_2+\|J\|_2)\sup_{j\in\Z}\int_0^t\|\Delta_jH(t')\|_\infty
dt'\big)\\
&&\qquad\qquad\qquad\qquad\qquad\qquad\qquad\quad\big(\sqrt{Ct}+\|H(0)\|_2+\|I(0)\|_2+\|J(0)\|_2\big).
\end{eqnarray}
Denote $$\zeta(T)\triangleq
\sup_{t\in[0,T)}\big(\|H(t)\|_2+\|I(t)\|_2+\|J(t)\|_2\big),$$ then
(3.24) can  be reduced to be
\begin{eqnarray}\label{3.24}
\zeta(T)\leq
\exp\bigg(C\log^+(\zeta(T))\sup_{j\in\Z}\int_0^t\|\Delta_jH(t')\|_\infty
dt'\bigg)\big(\sqrt{CT}+E(0)\big).
\end{eqnarray}
We should point out that the above inequality still holds if the
time interval $[0,T)$ is replaced by $[T-\varepsilon,T)$, that is
\begin{eqnarray}\label{3.26}
\zeta(T)\leq
\exp\bigg(C\log^+(\zeta(T))\sup_{j\in\Z}\int_{T-\varepsilon}^T\|\Delta_jH(t')\|_\infty
dt'\bigg) \big(\sqrt{C\varepsilon}+\zeta(T-\varepsilon)\big).
\end{eqnarray}
Setting $Z(T)\triangleq\log^+(\zeta(T))=\log(e+\zeta(T))$, thanks to
(\ref{3.26}), we have
\begin{eqnarray}\label{3.27}
Z(T)\leq \log\big(\sqrt{C\varepsilon}+\zeta(T-\varepsilon)+e\big)+
CZ(T)\sup_{j\in\Z}\int_{T-\varepsilon}^T\|\Delta_j(\nabla\times
u)(t')\|_\infty dt',
\end{eqnarray}
by the condition (\ref{1.2}) of Theorem\ref{the1},
$$\lim_{\varepsilon\rightarrow 0}\sup_{j\in \Z}
\int_{T-\varepsilon}^T\|{\Delta}_j(\nabla\times u)\|_\infty
dt=\delta<M,$$ we know that, when $\varepsilon\rightarrow 0$, if we
choose $MC$ is small enough, then it has
\begin{eqnarray}\label{3.28}
Z(T)\leq CZ(T-\varepsilon).
\end{eqnarray}
On the other hand, by multiplying $(u,\omega,b)$, it can be easily
derived  from Magneto-micropolar fluid equation(\ref{1.1}) that
\begin{eqnarray}\label{3.29}
\nonumber&&\|u\|_2^2+\|\omega\|_2^2+\|b\|_2^2+2\mu\int_0^t\|\nabla
u\|_2^2 dt'+2\gamma\int_0^t\|\nabla \omega\|_2^2
dt'+2\nu\int_0^t\|\nabla
b\|_2^2 dt'\\
&&+2\kappa\int_0^t\|div \omega\|_2^2 dt'+2\chi\int_0^t\|
\omega\|_2^2 dt'\leq \|u_0\|_2^2+\|\omega_0\|_2^2+\|b_0\|_2^2.
\end{eqnarray}
(\ref{3.29}) along with (\ref{3.28}) imply that
\begin{eqnarray}\label{3.30}
\nonumber&&\sup_{t\in[T-\varepsilon,T)}\big(\|u(t)\|_{H^1}+\|\omega(t)\|_{H^1}+\|b(t)\|_{H^1}\big)\\
&&\leq
C(\|u(T-\varepsilon)\|_{H^1}+\|\omega(T-\varepsilon)\|_{H^1}+\|b(T-\varepsilon)\|_{H^1}).
\end{eqnarray}
Hence by (\ref{3.30}) and (\ref{3.8}), we can get the $H^s$
regularity at time $t=T$, that is the  smooth  solution
$(u,\omega,b)$ can be extended past time $T$, that's the end of the
proof.$\square$

\begin{center}

\end{center}
\end{document}